\documentclass[a4paper]{article}
\usepackage{setspace}
\usepackage{graphicx}
\usepackage{epstopdf} 
\usepackage{epsfig}
\usepackage{amsmath}
\usepackage{amssymb} 
\usepackage{anysize}
\usepackage{algorithmic}
\usepackage{booktabs}
\usepackage{algorithm2e}
\usepackage{url}
\usepackage{rotating}
\usepackage{verbatim}

\usepackage{natbib}
\allowdisplaybreaks
\usepackage{pst-all}

\def\CC{{\mathcal{C}}}

\def\JJ{{\mathcal{J}}}
\def\RR{{\mathcal{R}}}

\def\SS{{\mathcal{S}}}
\def\FF{{\mathcal{F}}}

\def\LL{{\mathcal{L}}}

\def\keywords{\vspace{.5em}
{\textit{Keywords}:\,\relax%
}}

\title{Solution of minimum spanning forest problems with reliability constraints}

\author{\small Ida Kalateh Ahani$^{(1)}$, Majid Salari$^{(1)}$\footnote{Department of Industrial Engineering, Ferdowsi University of Mashhad, P.O. Box 91775-1111, Mashhad, Iran. TEL: +98 051 38805113, e-mail: msalari@um.ac.ir}, Seyed Mahmoud Hosseini$^{(1)}$, Manuel Iori$^{(2)}$\\
\small (1) Department of Industrial Engineering, Ferdowsi University of Mashhad, Mashhad, Iran\\
\small \url{{i_ahani, msalari, sm_hosseini}@um.ac.ir} \\
\small (2) Department of Sciences and Methods for Engineering, University of Modena and Reggio Emilia, Italy\\
\small \url{manuel.iori@unimore.it} \\
}

\begin{document}
\singlespacing
\date{}
\maketitle
\doublespacing

\begin{abstract}

We propose the reliability constrained $k$-rooted minimum spanning forest, a relevant optimization problem whose aim is to find a $k$-rooted minimum cost forest that connects given customers to a number of  supply vertices, in such a way that a minimum required reliability on each path between a customer and a supply vertex is satisfied and the cost is a minimum. 
The reliability of an edge is the probability that no failure occurs on that edge, whereas the reliability of a path is the product of the reliabilities of the edges in such path.
The problem has relevant applications in the design of networks, in fields such as telecommunications, electricity and transports. For its solution, we propose a mixed integer linear programming model, and an adaptive large neighborhood search metaheuristic which invokes several shaking and local search operators. Extensive computational tests prove that the metaheuristic can provide good quality solutions in very short computing times. 
\end{abstract}


\section{Introduction} \label{sec:Introduction} 

The \emph{minimum spanning tree} (MST) is one of the most celebrated problems in the field of combinatorial optimization. It is defined on a connected, undirected and weighted graph $G=(V,E)$, where \textit{V} is the set of vertices and \textit{E} is the set of edges. An MST of \textit{G} is a tree composed of a subset of edges of \textit{E} spanning all the vertices of $V$ and having minimum total weight. The MST has several applications, including design of computer and communication networks, picture processing, automatic speech recognition, clustering and classification problems (see, e.g., \citealt{WC04}). Several polynomial time algorithms have been proposed to solve it, starting from \cite{K56} and \cite{P57}.

A forest consists of a set of mutually disjoint trees, and a \emph{spanning forest} (SF) is a forest that spans all vertices of the graph. When the graph is connected, an SF is a forest where each vertex of the graph is included in one of the trees (see, e.g., \citealt{MV10, DSL10}). When instead the graph is not connected but consists of several connected components, an SF is a sub-graph containing a spanning tree of each component (see, e.g., \citealt{KB09, NCKB12}).
A \emph{\textit{k}-rooted spanning forest} is an SF with \textit{k} disjoint trees, where each tree has its own \emph{root}. These problems have different applications in designing supply chain networks, image processing (\citealt{TCB10, BTACB12}), political districting (\citealt{Y09}) and designing communication networks (\citealt{YTK96}), among others. 

\begin{figure}[htb]
  \caption{An illustrative example of the RCKRMSF problem.}
  \centering
    \includegraphics[width=80mm]{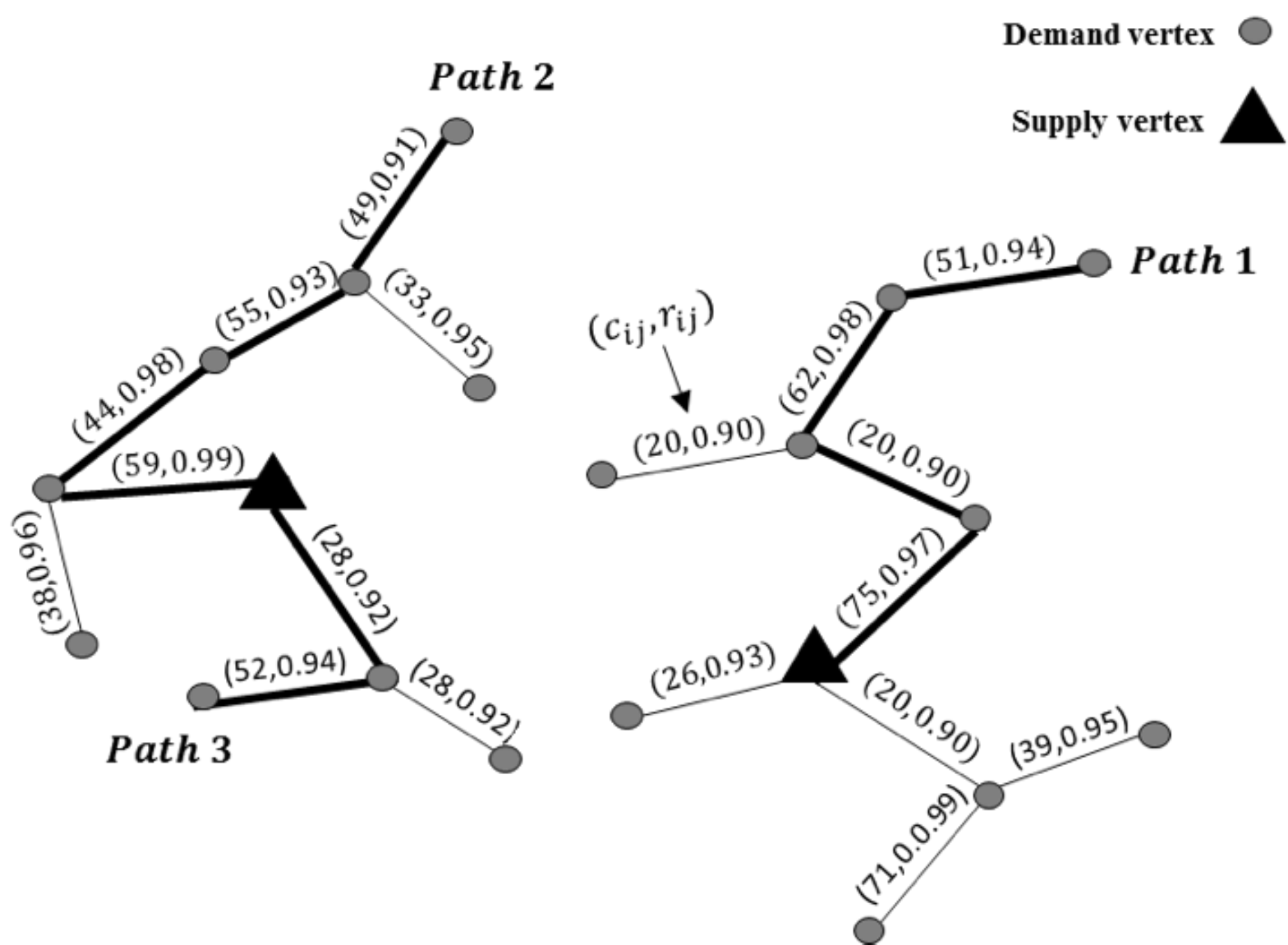}
    \label{fig:example}
\end{figure}

In this paper, we study  the \emph{reliability constrained \textit{k}-rooted minimum spanning forest} (RCKRMSF). Essentially, in the RCKRMSF we are given two sets of vertices, one including demand vertices and the other supply vertices. Demand vertices can be seen as customers and  should be connected by means of a path in the graph to one of the supply vertices. Each edge $(i, j)$ of the graph is associated with a cost $c_{ij}$ and a reliability $r_{ij}$.  The reliability of an edge is the probability that no failure occurs on that edge, and the reliability of a path is the product of the reliabilities of all edges included in the path. The goal of the RCKRMSF is to construct a \textit{k}-rooted minimum-cost SF that connects all customers to one or more supply vertices, and for which any path between a customer and its supply vertex satisfies a minimum required reliability.
{An RCKRMSF example is depicted in Figure \ref{fig:example}, where 2 supply and 18 demand vertices are rooted using an SF with 2 trees. Costs and reliabilities are reported close to each arc. Three paths are highlighted: Path 1 has cost 
75+20+62+51=208 and reliability 0.97$\times$0.90$\times$0.98$\times$0.94$\approx$0.80; similarly, path 2 has cost 207 and reliability 0.82 and path 3 cost 80 and reliability 0.86. It is easy to check that the solution guarantees minimum reliability 0.80 on all paths.}

The RCKRMSF is of interest because models many real-world applications in different fields. Reliability of paths is indeed a crucial component of systems arising in fields such as telecommunication, electricity and transport networks, just to cite some (see, e.g., \citealt{ENR15}). For example, when broadcasting information or data from servers to terminal nodes of a network, it is often needed to guarantee that the data is received  with at least a pre-specified reliability. Similarly, in radio communication networks, in case of a critical situation, the reliability of communication among central departments and stations should be greater than a minimum required level. Also, when a disaster occurs, reliability of delivering services to people in declared areas should be taken into account besides the minimization of delivery costs (\citealt{YSA18}). In all these examples, the design a network has to take into account both costs and minimum reliability, and could be modeled as an RCKRMSF. In general, problems where reliability matters in distributing services, data, information, or in general a given commodity from supply vertices to customers can be viewed as instances of the RCKRMSF.

Despite its large number of applications, to the best of our knowledge the RCKRMSF has not been previously addressed in the literature. In this paper, we first formally model the problem by using \emph{mixed integer linear programming} (MILP). As the MILP model is effective only on small-scale instances, we also propose a metaheuristic, based on the concept of \emph{adaptive large neighborhood search} (ALNS) and enriched with several local search routines and shaking operators. We prove the computational effectiveness of the proposed ALNS on a large set of instances. We also test the ALNS on other simplified versions of the RCKRMSF that have been tackled in the optimization literature, obtaining interesting insights. 

The paper is organized as follows. A concise literature review of the wide area of SF applications is provided in Section \ref{sec:review}. A formal problem description and a MILP model are presented in Section \ref{sec:Modeling}. The details of the proposed ALNS are given in Section \ref{sec:Solutionmethod}. Extensive computational tests on the RCKRMSF and other problem variants are discussed in Section \ref{sec:Results}, and conclusions are drawn in Section \ref{sec:Conclusion}.

\section{Brief literature review} \label{sec:review} 

We are not aware of previous works focused on the RCKRMSF, but we found a large number of works that address different constrained variations of the classical MST that are related to the RCKRMSF. 

The \emph{capacitated MST} (CMST) is an MST generalization in which each vertex $i$ of the graph, with the exception of a given supply vertex, has a non-negative demand $d_i$. The objective is to find an MST such that the sum of the demands in each sub-tree, rooted at the supply vertex, does not exceed a preset capacity value. The CMST is interesting because models many real-world situations and also because it can serve as a lower bound for capacitated vehicle routing (see, e.g., \citealt{TV14}) and pickup-and-delivery (see, e.g., \citealt{BI17}) problems. It has been tackled by many exact algorithms, including polyhedral approaches (\citealt{GL05}) and branch-and-cut-and-price algorithms (\citealt{UFLPPA08}). In terms of heuristics, good computational results were obtained with a biased random-key genetic algorithm by \cite{RAFR15} and with a quick greedy heuristic by \cite{KI17}.

\cite{NH80} proposed the \emph{degree constrained minimum spanning tree problem}, an MST variation in which a threshold on the maximum number of edges incident to each vertex is imposed. The problem models a number of applications, including the branching of cables in offshore wind farms (\citealt{KHBM15}). It has been addressed with several exact and heuristic techniques, such as variable neighborhood search (\citealt{SM08}), branch-and-cut (\citealt{CS14}) and branch-and-cut-and-price (\citealt{BCL16}). An analysis of problem variants as well as new MILP formulations have been recently presented by \cite{DCSA17}.

\cite{ACCG92} studied the \emph{bounded diameter MST} (BDMST). The objective of the BDMST is to find an MST such that the number of edges between any pair of vertices does not exceed a given input value. A relevant BDMST variation is the \emph{hop constrained MST} (HCMST), where each path between a designed root vertex and any other vertex should not contain more than $\delta$ edges (hops). Hop constraints can be seen as a first attempt to guarantee a certain level of service with respect to some performance constraints such as reliability. Indeed, if a certain reliability value $\beta$ is associated with each edge of the graph, then limiting the number of edges to $\delta$ has the effect of maintaining a minimum required level of reliability $\beta^{\delta}$ (see, e.g., \citealt{G95,L16}).

More elaborated constraints are adopted in the \emph{(rooted) delay constrained MST} (DCMST), also known as  the \emph{(rooted) distance constrained MST}. The DCMST is a generalization of the HCMST in which a delay $l_{ij}$ is associated with each edge $\left\{i, j\right\} \in E$, and the sum of all delays in a path should not be greater than a pre-specified threshold $L$ (see, e.g., \citealt{SRV97, GPS08}). 

Another aspect that is relevant for the RCKRMSF is that the optimized solution is divided into clusters, at most one for each supply vertex. Clustering has been intensively studied in the area of network design. Here we only mention the related \emph{generalized MST} (GMST), originally proposed by \cite{MLT95}. The GMST calls for finding an MST that either contains exactly one vertex for each cluster (as in, e.g., \citealt{MLT95,PMSP18}), or alternatively contains at least one vertex for each cluster (as in \citealt{DHC00}). 

As previously mentioned, reliability is the probability of failure-free performance of a system for a specified duration and under stated conditions. 
There are many measurements for the reliability of a network, such as two-terminal, $k$-terminal and all-terminal reliability. Two-terminal reliability is defined as the reliability between two specified vertices, whereas the reliability of communication between every pair of vertices in the network is called all-terminal reliability (see, e.g., \citealt{KSI01, SKS02}). 
A general case, covering the two previous cases, is the $k$-terminal reliability, in which minimum reliability in the designed network must be respected among all pairs of a specified subset of $k$ vertices (see, e.g., \citealt{CCRRS13, RDD16, HSS18}).

In designing communication networks, both reliability and installation costs are of great importance, so many researchers have taken into account both these aspects when developing algorithms or models to design such networks. In general, designing networks with the presence of budget and reliability constraints can be achieved in two different ways: the objective can be the cost minimization subject to minimum reliability guarantee, or the reliability maximization subject to a maximum budget limit.
In our work, we consider cost minimization subject to minimum reliability. For the literature on this type of problems, we mention the branch-and-bound algorithm by \cite{JHC93},  the method based on  cross entropy by \cite{AD09}, and the ant colony optimization approaches by \cite{WW09} and \cite{DAB10}.
For what concerns instead reliability maximization, we mention
the simulated annealing by \cite{AR93}, the hill-climbing, simulated annealing and genetic algorithms by \cite{ADS03}, the so-called shrinking and searching algorithm by \cite{SSH05}, and the artificial neural network by \cite{DBTT12}.
Bi-objective models have been also studied to take care of the cost and reliability components in an integrated manner. Among these, we cite the genetic algorithms by \cite{APKKS09} and \cite{KS11}. We also mention the recent label setting algorithm by \cite{SPG18}, which deals with multiple objectives and could be adapted to take care of reliability. 

We refer readers interested in further pointers to the literature to the surveys by \cite{BCP95} on network reliability, \cite{BH01} and \cite{AK11} on MST algorithms, and \cite{AFC15} on multi-objective models for reliability.


\section{Problem formulation and mathematical model} \label{sec:Modeling} 

In this section, we provide a formal RCKRMSF description and a MILP model. The model is based on a multicommodity network flow (see, e.g., \citealt{MW95}), where each commodity is assigned to a given customer and its flow is used to model the path connecting the customer to one of the supply vertices. 
To this aim, we introduce a directed graph $G'=(V,A)$, obtained by replacing each edge $\left\{i,j\right\} \in G$ with two directed arcs $\left(i, j\right)$ and $\left(j,i\right)$. The vertex set is partitioned as $V = S \cup D$ and the arc set is defined as $A = \left\{(i,j): i,j \in V\right\}$. Sets $S = \left\{1,2,\dots, m\right\}$ and $D = \left\{m+1, m+2,\dots,n\right\}$ represent the supply and demand vertices, respectively. Two positive values, namely the cost $c_{ij}$ and the reliability $r_{ij}$, are associated with each arc $(i,j) \in A$. 
Recall that the reliability of a path is the product of the reliabilities of all arcs in the path, and let $\alpha$ be a minimum required level of reliability.
The objective of the RCKRMSF is to construct a $k$-rooted minimum-cost SF of the vertices in \textit{V}, such that each demand vertex (also denoted as customer) is connected by a path in the SF to exactly one of the supply vertices, and by ensuring that each such path has reliability at least equal to $\alpha$.  Without loss of generality, to ensure that a feasible solution exists we suppose that each demand vertex $i \in D$ can be linked to at least one source vertex $j \in S$ by a direct arc having reliability $r_{ij} \geq \alpha$.

The MILP model that we introduce makes use of two families of binary variables and a family of continuous non-negative variables, namely:
\begin{eqnarray}
& y_{i}^{s} =  \left\{ \begin{array}{ll}
1 \quad \mbox{if customer $i$ is allocated to supply vertex $s$},\enspace\\
0 \quad \mbox{otherwise}\end{array}\right. & \forall i \in D, s \in S, \label{eq:var1}\nonumber \\
& x_{ij}^{s} = \left\{ \begin{array}{ll}
1 \quad \mbox{if arc $(i,j)$ belongs to a path originated from supply vertex $s$},\enspace\\
0 \quad \mbox{otherwise}\end{array}\right. & \forall (i,j) \in A, s \in S, \label{eq:var2}\nonumber \\
& f_{ij}^{h} = \mbox{flow of commodity $h$ along arc $(i,j)$} & \forall (i,j) \in A, h \in D. \nonumber
\end{eqnarray}

A commodity is assigned to each customer $h \in D$, and the flow variables $f$ are used to build paths that start at one of the supply vertices and end at the customers. In particular, a unit of flow of each commodity $h$ is imposed to leave a supply vertex and reach $h$. Our model is then as follows:
\begin{eqnarray}
\textrm{(RCKRMSF)}: \label{eq:obj1} \min \sum_{s\in S}\sum_{\left(i, j\right)\in A} c_{ij}x_{ij}^{s}   \enspace\\
\label{eq:con1} \sum_{s\in S}\sum_{\left(i, j\right)\in A} x_{ij}^{s} = n - m    \,\,\,                                         				 \enspace\\
\label{eq:con2} \sum_{s\in S} (x_{ij}^{s} + x_{ji}^{s})               \le 1      \,\,\, &&  \forall \left\{i, j\right\}\in E, 							 \enspace\\
\label{eq:con3} \sum_{s\in S} y_{i}^{s}                               =   1      \,\,\, &&  \forall i\in D,                							 \enspace\\
\label{eq:con4}                  y_{s}^{s}                               =   1   \,\,\, &&  \forall s\in S,                							 \enspace\\
\label{eq:con5}                  x_{ij}^{s}                         \le y_{i}^{s}\,\,\, &&  \forall \left(i,j\right)\in A: i\in D, s\in S,\enspace\\
\label{eq:con6}                  x_{ij}^{s}                         \le y_{j}^{s}\,\,\, &&  \forall \left(i,j\right)\in A: j\in D, s\in S,\enspace\\
\label{eq:con7} \sum_{j\in D}f_{sj}^{h}                        = y_{h}^{s}  \,\,\, &&  \forall s\in S, h\in H,        		 					 \enspace\\
\label{eq:con8} \sum_{i\in V}f_{ih}^{h} - \sum_{j\in D}f_{hj}^{h} = 1\,\,\, &&  \forall h\in H,                		 					 \enspace\\
\label{eq:con9} \sum_{i\in V}f_{ij}^{h} - \sum_{i\in D}f_{ji}^{h}  = 0      \,\,\, &&  \forall h\in H, j\in D, j\neq h,        		  		 \enspace\\
\label{eq:con10} f_{ij}^{h} \le \sum_{s\in S}x_{ij}^{s}             \,\,\, &&  \forall h\in H, \left(i,j\right)\in A,      					 \enspace\\
\label{eq:con11} \prod_{\left(i, j\right)\in A} r_{ij}^{f_{ij}^{h}} \geq \alpha  \,\,\, &&  \forall t_h\in D,															 \enspace\\
\label{eq:con12} x_{ij}^{s} \in \{0,1\}  																				 \,\,\, &&  \forall \left(i,j\right)\in A, s\in S, 	 		 \enspace\\
\label{eq:con13} y_{i}^{s}  \in \{0,1\} 																				 \,\,\, &&  \forall i\in D, s\in S,                	 		 \enspace\\
\label{eq:con14} f_{ij}^{h} \ge 0 																							 \,\,\, && \forall \left(i,j\right)\in A, h\in H.  	 		 \enspace
\end{eqnarray}

The objective function \eqref{eq:obj1} minimizes the total cost of the SF. Constraint \eqref{eq:con1} forces the SF to contain exactly $n-m$ arcs. Constraints \eqref{eq:con2} guarantee that each edge $\left\{i,j\right\}\in E$ is assigned to at most one supply vertex and  carries flow in at most one direction. Constraints \eqref{eq:con3} impose that each customer is assigned to exactly one supply vertex. Constraints \eqref{eq:con4} simply require that a supply vertex is allocated to itself. Constraints \eqref{eq:con5} and \eqref{eq:con6} ensure that if arc $(i, j)$ is assigned to a path originating from $s\in S$, then both vertices \textit{i} and \textit{j} are allocated to $s$. 
Constraints \eqref{eq:con7}--\eqref{eq:con9} are used to impose flow conservation, in particular: the flow of commodity $h$ leaving supply vertex $s$ is set to be equal to $y_h^s$; 
the difference between the amount of commodity $j$ entering and leaving a customer $h$ is set to be 1 when $j=h$, and 0 when $j \neq h$.
Constraints \eqref{eq:con10} link together $x$ and $y$ variables, by ensures that if a flow is sent along an arc, then such arc must be included in the solution. 
Constraints \eqref{eq:con11} impose the minimum reliability level $\alpha$ on each path in the solution. Then, constraints \eqref{eq:con12}--\eqref{eq:con14} define the domains of the variables.

Constraints \eqref{eq:con11} are non-linear, but can be linearized by using the logarithmic function (as in, e.g., \citealt{WA88}). Essentially, we can replace the left-hand side and the right-hand side of the constraint by their logarithms, obtaining:
\begin{eqnarray}
\label{eq:con15}\sum_{\left(i, j\right)\in A} f_{ij}^{h}\ln\left(r_{ij}\right) \geq \ln \left(\alpha\right) \,\,\, &  \forall h\in D.\enspace
\end{eqnarray}

The resulting MILP model, composed of \eqref{eq:obj1}--\eqref{eq:con10} and \eqref{eq:con12}--\eqref{eq:con15}, has been solved by means of a commercial solver, and the results that have been obtained are provided in Section \ref{sec:Results} below.

\section{Adaptive Large Neighborhood Search} \label{sec:Solutionmethod} 


The {\em large neighborhood search} (LNS) algorithm, originally proposed by \cite{S98}, attempts at improving a given solution by alternately destroying and repairing it (\cite{MS18}). 
The ALNS is an extension of the LNS that has been first proposed by \cite{RP06} and since then has obtained successful results for a large variety of optimization problems (\citealt{PSL17, AS18}). 
The ALNS allows the use of multiple destroy and repair algorithms within the same search, and adopts an adaptive selection mechanism to decide which operators should be used. Essentially, at each iteration a given part of the solution is destroyed and repaired by applying the appropriate destroy and repair rules, and these rules are selected on the basis of given probabilities that depend on their performance during the search process. 

In this section, we propose a heuristic algorithm for the RCKRMSF that is based on the ALNS paradigm and follows the pseudocode given in Algorithm \ref{sec:Alg1}.
We are given a set of $n_1$ (resp. $n_2$) local search (resp. shaking) procedures. A weight $w^{k}_{ls}$ (resp. $w^{k}_{sh}$) is associated with each local search (resp. shaking) procedure $k$, and is initially set to 1 for all procedures. With each local search (resp. shaking) procedure $k$, we also associate a score $\pi^{k}_{ls}$ (resp. $\pi^{k}_{sh}$) representing the performance of the procedure during the search process. 

Our ALNS contains a given number $\phi_1$ of global iterations, called segments for short, each consisting of $\phi_2$ inner iterations. 
At the beginning of each segment, all scores are set to zero. Then, at each iteration, a local search procedure $i$ and a shaking procedure $j$ are chosen by the roulette-wheel mechanism with  probability $p_i = w^{i}_{ls}/\sum_{k=1}^{n_1}w^{k}_{ls}$ and $ p_j = w^{j}_{sh}/\sum_{k=1}^{n_2} w^{k}_{sh}$, respectively. 
A new solution is obtained by applying first the shaking and then the local search on a current solution. Then, the scores are updated as follows: if the cost of the new solution is better than that of the incumbent, the scores of the involved local search and shaking procedures are augmented by $\zeta_1$; if the cost of the new solution is better than that of the previous temporary solution but not better than that of the incumbent, the scores are augmented by $\zeta_2$;  if the new solution is non-improving but accepted by a {\em Simulated Annealing} (SA) criterion, the scores are increased by $\zeta_3$. 
Let $Cost(X)$ be the cost of solution \textit{X}. The SA criterion works as follows: given a current solution $CurrentSolution$, a neighbor solution $TempSolution$ is always accepted if $Cost(TempSolution) < Cost(CurrentSolution)$; otherwise, it is accepted with probability $e^{(Cost(TempSolution) - Cost(CurrentSolution))/\theta}$, in which $\theta$ is the current temperature. The initial temperature is set to $\theta_{0}$ and is decreased gradually, by performing the update $\theta_{new} = \lambda\theta_{old}$ where $0<\lambda<1$ is the cooling rate. In addition, the final temperature is set to $\epsilon$. 

When a segment ends, the weights associated with the local search and shaking procedures are updated according to the scores obtained during the last $\phi_2$ iterations using the following equation: 
\begin{equation}
w^{k}_{\mu}=\left\{ \begin{array}{ll}
w^{k}_{\mu} & \quad \mbox{if $\gamma^{k}_{\mu}$ = 0},\enspace\\
(1-\eta)w^{k}_{\mu}+ \eta\pi^{k}_{\mu}/\gamma^{k}_{\mu}  & \quad \mbox{ otherwise.
}\end{array}\right. \quad \mu = sh, ls \label{eq:var5} \end{equation}
In \eqref{eq:var5}, $\eta \in \left[0,1\right]$ is an input parameter called reaction factor, whereas $\gamma^{k}_{ls}$ (resp. $\gamma^{k}_{sh}$) gives the number of times a local search (resp. shaking) procedure $k$ has been selected during the $\phi_2$ inner iterations of the last segment.

\begin{table}
\begin{algorithm}[H]
\hrulefill \\
\small \label{sec:Alg1}
\textbf{Inputs}:
$\theta$: initial temperature; $\epsilon$: final temperature; $\lambda$: cooling rate; $\eta$: reaction factor; $n_1$: number of local search procedures; $n_2$: number of shaking procedures\;
\textbf{Output}: \textit{BestSolution}\;
\textit{TempSolution}     =  \textit{CurrentSolution} =  \textit{BestSolution}   = \textit{Initialization}()\;
\textbf{for each} \textit{local search} $l\in \left\{1,2,\cdots,n_1\right\}$ set $w^{l}_{ls}$ = 1\;
\textbf{for each} \textit{shaking procedure} $k\in \left\{1,2,\cdots,n_2\right\}$ set $w^{k}_{sh}$ = 1\;
\For{$($$Segment_{iter} = 1$ \textbf{to} $\phi_1$$)$}{
  \textbf{for each} \textit{local search} $l\in \left\{1,2,\cdots,n_1\right\}$ set $\pi^{l}_{ls} = 0$ \textbf{and} $\gamma^{l}_{ls}$ = 0\;
  \textbf{for each} \textit{shaking procedure} $k\in \left\{1,2,\cdots,n_2\right\}$ set $\pi^{k}_{sh} = 0$ \textbf{and} $\gamma^{k}_{sh}$ = 0\;
 	\For{$($$Local_{iter} = 1$ \textbf{to} $\phi_2$$)$}{
 	      \textit{Local search} =  choose a procedure \textit{i} using roulette wheel with weights $w^{l}_{ls}$ ($l\in \left\{1,2,\cdots,n_1\right\}$)\;
 	      \textit{Shaking} = choose a {procedure} \textit{j} using roulette wheel with weights $w^{k}_{sh}$ ($k\in \left\{1,2,\cdots,n_2\right\}$)\;
 	      $\gamma^{i}_{ls}$ += 1; $\gamma^{j}_{sh}$ += 1\;
 	      \textit{TempSolution} =  \textit{Shaking}(\textit{TempSolution})\;
 	      \While{$($TempSolution \emph{can be improved}$)$}{
 	            \textit{TempSolution} =  \textit{Local search}(\textit{TempSolution})\;
 	      }
 	      \If {$($Cost$($TempSolution$)$ $<$ Cost$($BestSolution$)$$)$} {
  					\textit{CurrentSolution} = \textit{BestSolution}    = \textit{TempSolution}\;
    				$\pi^{i}_{ls}$ += $\zeta_{1}$; $\pi^{j}_{sh}$ += $\zeta_{1}$\;
  			} \ElseIf{$($Cost$($TempSolution$)$ $<$ Cost$($CurrentSolution$)$$)$} {
  							\textit{CurrentSolution}  =  \textit{TempSolution}\;
  							$\pi^{i}_{ls}$ += $\zeta_{2}$; $\pi^{j}_{sh}$ += $\zeta_{2}$\;
  			}
  			\ElseIf {$($$\theta$ $>$ $\epsilon$$)$} {
    			    	  $\rho$ = $e^{(Cost(BestSolution)-Cost(TempSolution))/\theta}$\;
    			    	  $\tau$ = rand (0,1)\;
    			    	  \eIf {$($$\tau$ $>$ $\rho$$)$} {
    			    	  		\textit{TempSolution} = \textit{BestSolution}\;
    			    	  }
    			    	   {
    			    	  			\textit{CurrentSolution} = \textit{TempSolution}\;
    			    	  			$\pi^{i}_{ls}$ += $\zeta_{3}$; $\pi^{j}_{sh}$ += $\zeta_{3}$\;
    			    	  }
									$\theta$ = $\theta \times \lambda$\;
    			}
    			\Else {
    						\textit{TempSolution} = \textit{BestSolution}\;
    		}
 	}
 	\textbf{for each} \textit{local search procedure} $l\in \left\{1,2,\cdots,n_1\right\}$ update $w^{l}_{ls}$ using \eqref{eq:var5}\;
 	\textbf{for each} \textit{shaking procedure} $k\in \left\{1,2,\cdots,n_2\right\}$ update $w^{k}_{sh}$ using \eqref{eq:var5}\;
}
\hrulefill
\caption{\small ALNS procedure for the RCKRMSF}
\end{algorithm}
\end{table}

In what follows, we provide the details of the inner procedures invoked by Algorithm \ref{sec:Alg1}. 

\subsection{Initialization} \label{sec:Initialization} 
To construct a feasible initial solution, we adopted an initialization procedure that applies the Prim algorithm over a modified version of the input graph, namely: we set to 0 all costs associated with arcs connecting two supply vertices; we execute the Prim algorithm on the resulting graph, so obtaining an MST; we remove all arcs connecting two supply vertices from the solution, so obtaining an SF. Then, if the SF is infeasible because does not satisfy the minimum reliability constraint,  the procedure proceeds to achieve feasibility as follows. Let $\LL$ represent the set of all leaf vertices of the infeasible SF. The procedure randomly selects a vertex $l \in \LL$ for which the reliability of the path (i.e., the path  from the selected supply vertex to \textit{l}) is not satisfied. Vertex \textit{l} is removed from its position in the current SF and allocated to a new position that is feasible and has minimum reinsertion cost. Note that this position always exists because we assumed in Section \ref{sec:Modeling} that there is at least a direct arc connecting a source vertex to a demand vertex having reliability greater than or equal to $\alpha$. The procedure is reiterated until all vertices are assigned to feasible positions, thus producing a feasible solution in output.

\subsection{Local search procedures} \label{sec:RepairAlgorithms} 
\begin{figure}[h!]
  \caption{An illustrative example of \textit{local search} 1.}
  \centering
    \includegraphics[angle=270,width=110mm]{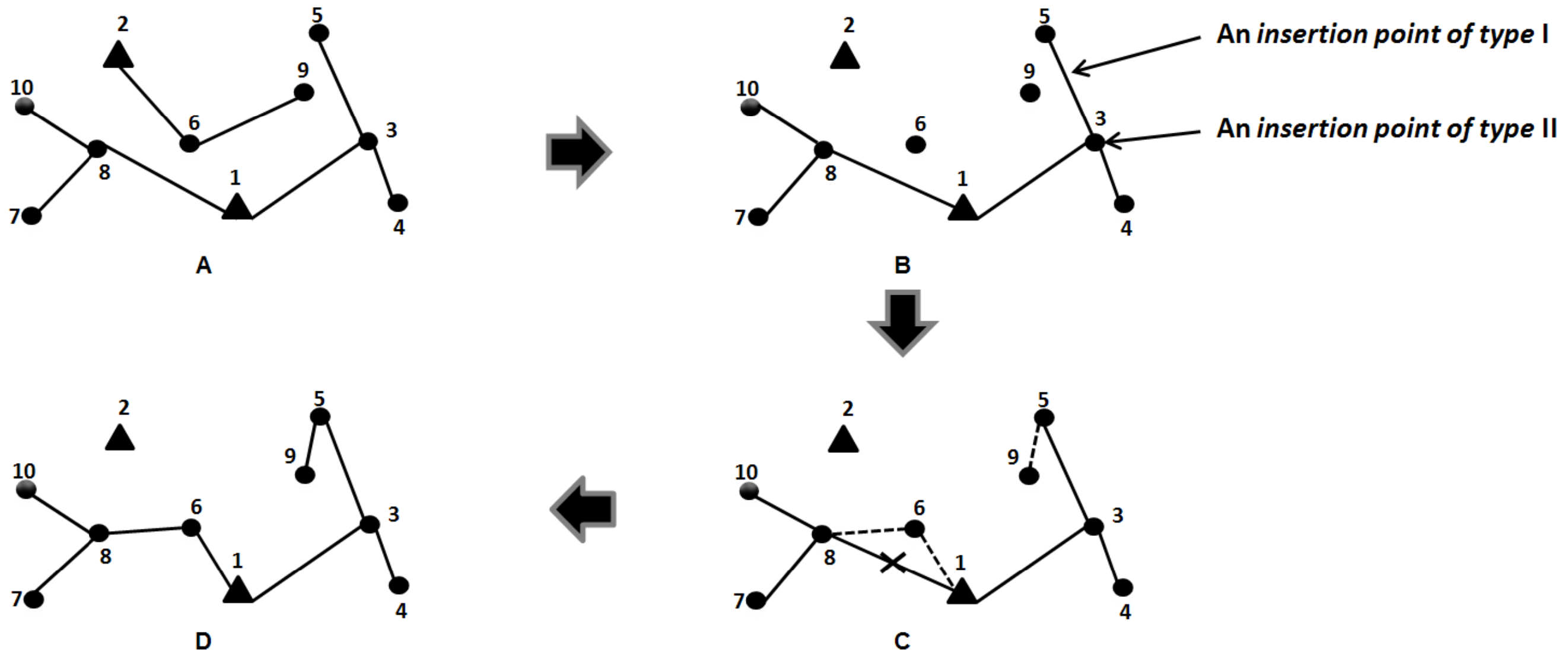}
    \label{fig:Repair1}
\end{figure}

Given an initial feasible solution $\FF$, five different local search operators have been implemented to try to improve the solution cost. Essentially, in all such procedures modifications of different parts of the initial solution are examined to possibly find a better solution by reconnecting demand and supply vertices.
\begin{itemize}
\item \textit{Local search} 1: a restricted solution $\RR\SS$ is obtained by removing an arc, chosen in random order, from the initial solution $\FF$. As a result, a subset $\LL^{'}$ of  vertices are disconnected from the supply vertices. The procedure examines relocating each vertex in $\LL^{'}$, one at a time, to feasible insertion points in $\RR\SS$. Essentially, insertion points can be an arc of $\RR\SS$ (\textit{insertion type} I) or a vertex (\textit{insertion type} II). A vertex $l\in \LL^{'}$ is selected in random order and allocated to its lowest-cost feasible insertion point in $\RR\SS$. The procedure iterates until all vertices in $\LL^{'}$ are relocated. In case the new solution has a better objective value than $\FF$, it is accepted, otherwise the procedure proceeds by selecting a new random arc of $\FF$. The process is reiterated until all arcs have been examined and no further improvement is possible. An example is represented in Figure \ref{fig:Repair1}:  Vertices 6 and 9 are disconnected from the solution by removing arc $(2,6)$ (Figures \ref{fig:Repair1}-A and \ref{fig:Repair1}-B) and allocated to new insertion points of types II and I, respectively (Figures \ref{fig:Repair1}-C and \ref{fig:Repair1}-D).
\end{itemize}

\begin{itemize}
\item \textit{Local search} 2: the procedure chooses a random edge and removes it from $\FF$. As a result, a sub-tree is disconnected in the restricted solution $\RR\SS$ that has been obtained. All insertion points of type II in $\RR\SS$ are examined for the re-insertion of the entire sub-tree. If a feasible and lower cost insertion point is found, then the sub-tree is moved there, otherwise it is put back in its original position. The procedure iterates until all arcs in $\FF$ have been tested for removal.
\end{itemize}

\begin{itemize}
\item \textit{Local search} 3: let $\LL$ represent the set of all leaves of  $\FF$. 
The procedure selects a random leaf from $\LL$ and removes it from its current position. It then attempts to re-insert it by considering all feasible insertion points of type I or II, and selecting the one having minimum cost and being feasible, if any. The move is accepted if it leads to an improvement in the objective function of the problem, and the procedure is re-iterated until all leaves have been examined and no further improvement is possible.
\end{itemize}

\begin{itemize}
\item \textit{Local search} 4: the procedure randomly selects two vertices of $\FF$ and swaps them. The move is accepted if the resulting solution is feasible and has a lower objective value. The procedure terminates when all possible swaps are examined and no further improvement is possible.
\end{itemize}

\begin{itemize}
\item \textit{Local search} 5: This rule is a math-based procedure that aims at improving the objective function by finding a new rearrangement of the leaves' positions in the initial solution $\FF$. The procedure starts by selecting all leaves of $\FF$ and removing them from the solution. A restricted solution $\RR\SS$ is consequently obtained. Assuming $\LL$ to represent the set of all extracted vertices, the procedure tries to find an improved cost solution by reallocating the vertices of $\LL$ to $\RR\SS$ by means of an integer linear programming model. 

To this aim, all sequences consisting of one or two customers belonging to $\LL$ are generated. Let $\SS'$ and $\SS''$ denote the sequences of size one and two of the extracted vertices, respectively. We represent a sequence $s \in \SS = \SS' \cup \SS''$ as $s = \left(b_s, t_s\right)$, in which $b_s$ and $t_s$ are the vertices belonging to \textit{s}, and $b_s = t_s$ in case $s\in S'$. The set of all feasible insertion points corresponding to the sequence $s$ is represented by $\JJ_s$. Essentially, an insertion point is a demand or a supply vertex in $\RR\SS$ to which sequence \textit{s} could be allocated through an \textit{insertion type} II move. We denote $re_j$ the reliability of the path ending at vertex \textit{j} and starting from one of the available supply vertices. In addition, we denote $re_s$ the reliability of the sequence $s \in S$. If $s\in \SS'$ the corresponding reliability is set to one, otherwise it is set to the reliability of arc $(b_s, t_s)$. Finally, we represent by $\SS_l$  the subset of sequences in $\SS$ that contain customer vertex $l \in \LL$.
We make use of the following decision variable:
\begin{equation}
y_{sj}=\left\{ \begin{array}{ll}
1 & \quad \mbox{ if sequence $s\in \SS$ is connected to the insertion point $j\in \JJ_s$},\enspace\\
0 & \quad \mbox{ otherwise,
}\end{array}\right.\label{eq:var3}\nonumber \end{equation}
for all $s \in \SS, j \in \JJ_s$. The model for the optimal reallocation of the leaves is then:
\begin{eqnarray}
\label{eq:obj2} \min \sum_{s\in S}\sum_{j\in \JJ_s} c_{sj}y_{sj}   \enspace\\
\label{eq:con18} \sum_{s\in \SS_l}\sum_{j\in \JJ_s} y_{sj} = 1   \,\,\, &  \forall l\in \LL,            \enspace\\
\label{eq:con19} re_j(1-y_{sj}) + re_j\times re_{s}\times r_{jb_s}\times y_{sj}   \geq \alpha              \,\,\, &  \forall j\in \RR\SS, s\in \SS, 	\enspace\\
\label{eq:con20} y_{sj} \in \{0,1\}  													 \,\,\, &  \forall  s\in \SS, j\in \JJ_s.
\end{eqnarray}
The objective function \eqref{eq:obj2} is to minimize the total reallocation cost. For each $l\in \LL$, constraints \eqref{eq:con18} impose that each extracted vertex is allocated to exactly one position of the restricted solution. Constraints \eqref{eq:con19} impose the minimum required reliability. Essentially, the reliability of any path obtained by adding a sequence $s\in \SS$ to a vertex $j\in \RR\SS$ cannot be less than the pre-specified value $\alpha$. Constraints \eqref{eq:con20} define the variables' domain.
\end{itemize}

\subsection{Shaking procedures} \label{sec:DestroyAlgorithms} 
In this section, we describe six shaking procedures that are used within the developed ALNS. In all such procedures, regardless of the change in the cost of the solution, new ways of connecting the vertices are considered, so that a larger portion of the feasible space can be searched and new solutions can be given in input to the local search operators above. All procedures start from an initial feasible solution $\FF$.

\begin{itemize}
	\item \textit{Shaking procedure} 1: The purpose of this procedure is to change a portion of $\FF$ by substituting an available path with a new one. To this aim, a demand vertex is randomly selected and the {\em max-reliability} path is obtained. Essentially, the max-reliability path is the path that originates from one of the supply vertices, reaches the selected demand vertex and has maximum reliability among all paths of this type. This path can be easily obtained by applying the Dijkstra algorithm. In case of having several paths with the same maximum reliability, one of them is randomly selected. The selected path is included in the solution, and the cycles that have been possibly generated by its addition are erased by removing arcs from $\FF$. The new solution obtained by this shaking procedure is always feasible from the reliability point of view. Figure \ref{fig:Destroy1} gives an example in which vertex 7 is selected (Figure  \ref{fig:Destroy1}.A) and its corresponding max-reliability path is calculated (dotted path in Figure  \ref{fig:Destroy1}.B). As the newly generated solution has cycles, the procedure removes a subset of arcs belonging to the arc set of the solution before adding the max-reliability path and connected to the vertices visited by the new path (Figure \ref{fig:Destroy1}.C), consequently restoring feasibility (Figure \ref{fig:Destroy1}.D).
\begin{figure}[h!]
  \caption{An illustrative example of \textit{shaking procedure} 1.}
  \centering
    \includegraphics[angle=270,width=110mm]{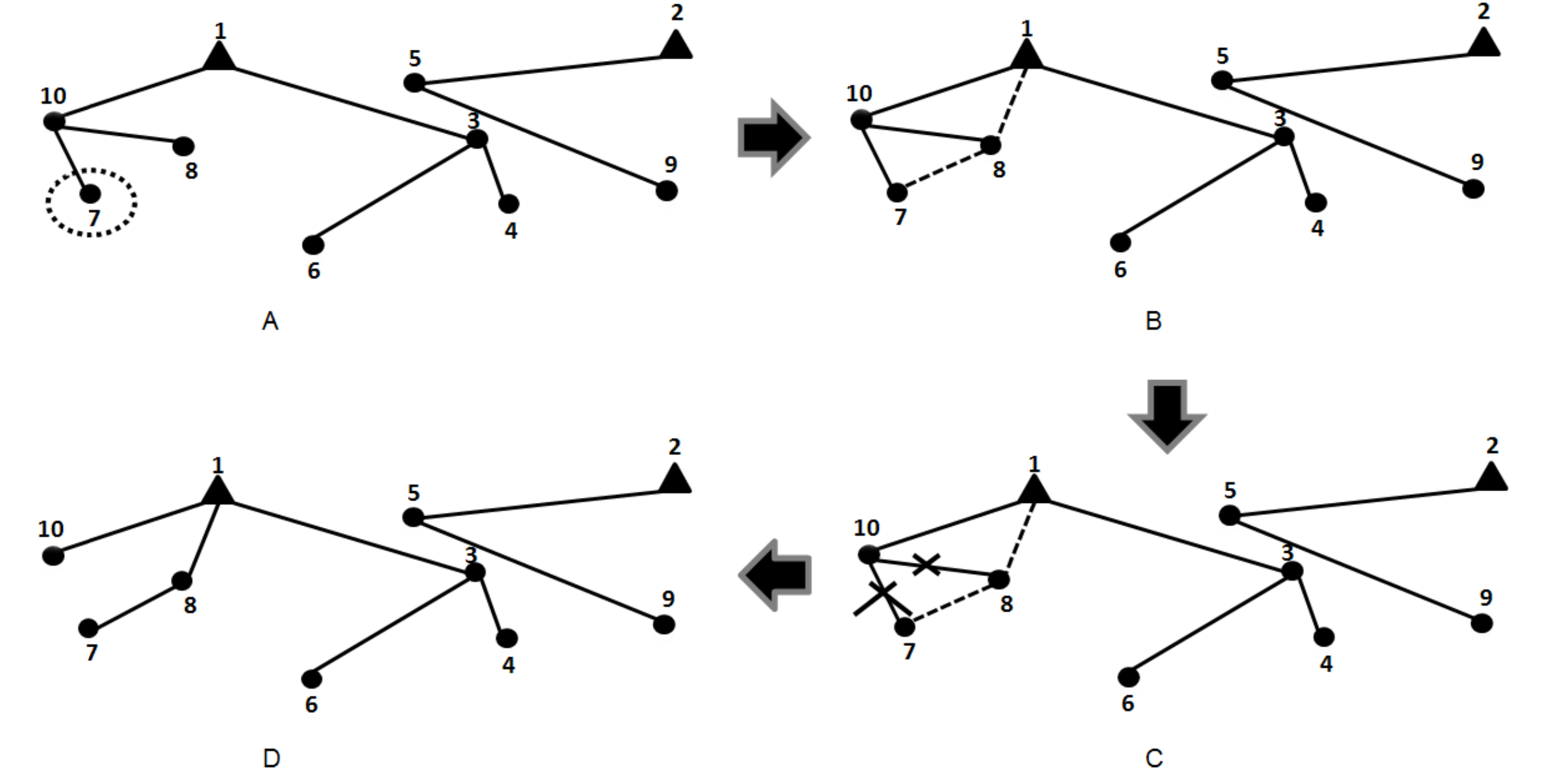} 
    \label{fig:Destroy1}
\end{figure}
\end{itemize}

\begin{itemize}
	\item \textit{Shaking procedure} 2: This procedure is a generalization of \textit{shaking procedure} 1, in which three demand vertices, instead of one, are randomly selected from $\FF$. Taking into account all the available supply vertices, the max-reliability paths for the three vertices are obtained and implemented in the new solution,  and possibly generated cycles are removed.
\end{itemize}

\begin{itemize}
	\item \textit{Shaking procedure} 3: We select two demand vertices, say \textit{i} and \textit{j}, having distance ($c_{ij}$) greater than a pre-specified input value $\CC$. The maximum-reliability paths for the two vertices are obtained and implemented, and possible cycles are removed from the new solution that has been obtained. An illustrative example is shown in Figure \ref{fig:Destroy3}, in which demand vertices 4 and 23 are selected to perform the shaking.
\begin{figure}[h!]
  \caption{An illustrative example of \textit{shaking procedure} 3.}
  \centering
    \includegraphics[angle=270,width=150mm]{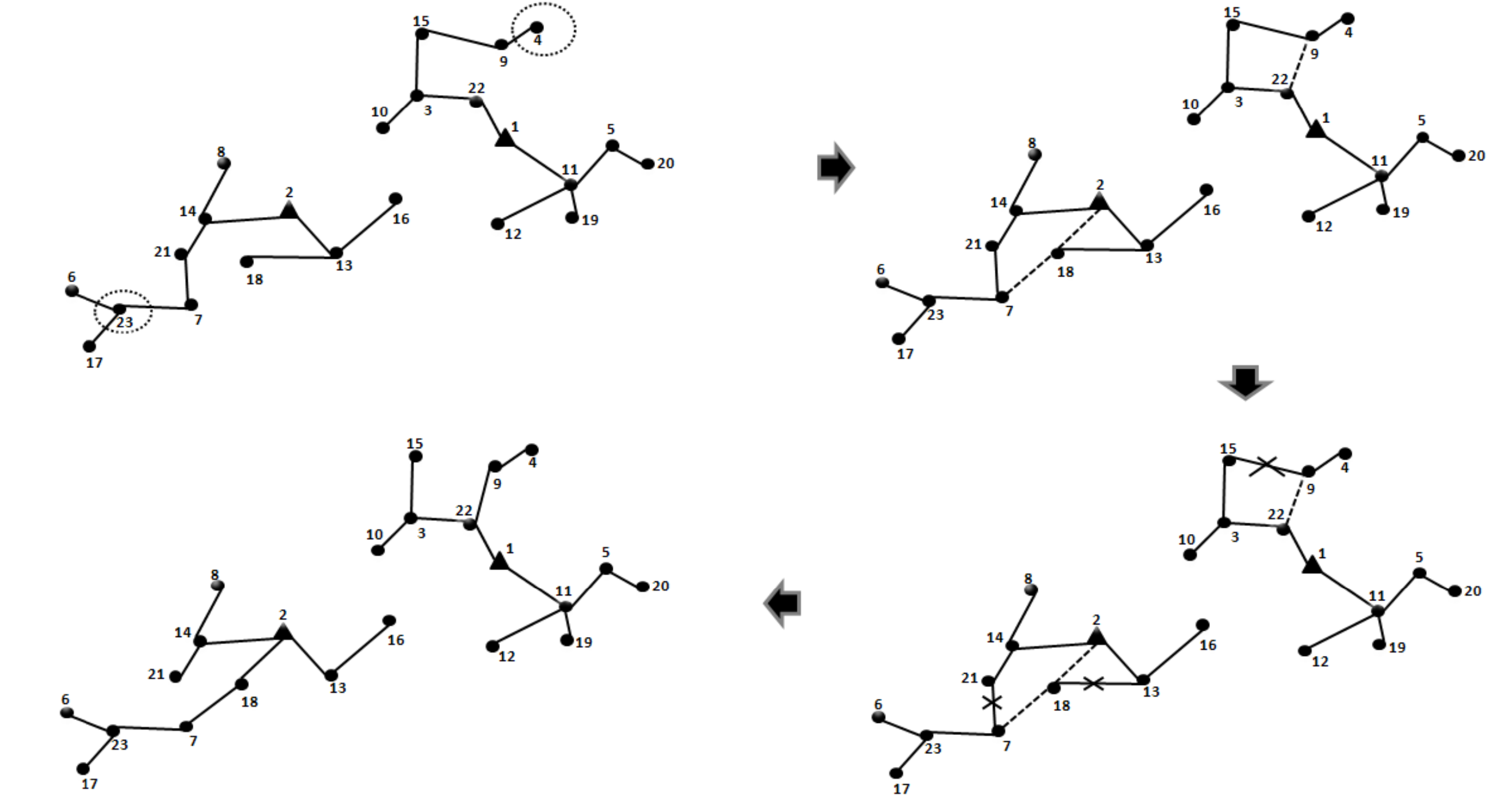}
    \label{fig:Destroy3}
\end{figure}
\end{itemize}

\begin{itemize}
\item \textit{Shaking procedure} 4: In this procedure, two vertices of $\FF$ are selected randomly, and their locations are swapped with each other. If the move is feasible from the reliability point of view, it is accepted, otherwise the procedure examines swapping other vertices to find the first feasible swap, if any.  
\end{itemize}

\begin{itemize}
\item \textit{Shaking procedure} 5: A restricted solution $\RR\SS$ is obtained by randomly removing an arc from $\FF$. As a result, a sub-tree $\FF^{'}$ of $\FF$ is disconnected from the forest. Then, a random vertex of $\RR\SS$ is selected to reallocate the extracted sub-tree $\FF^{'}$. If the solution obtained in this way is not feasible because does not satisfy the minimum reliability constraint, another random vertex of $\RR\SS$ is selected for the reallocation of $\FF^{'}$. In case no feasible position exists to reinsert $\FF^{'}$, $\FF^{'}$ is reallocated back to its original position, and the procedure iterates by selecting and removing a new random arc from $\FF$.
\end{itemize}
\begin{itemize}
\item \textit{Shaking procedure} 6: Similarly to what proposed for {shaking procedure} 5, in this procedure too an arc is randomly chosen from $\FF$ and removed. The disconnected component $\FF^{'}$ obtained by this operation is then reattached to $\RR\SS$. The only difference with respect to the previous procedure is that now the position in which $\FF^{'}$ has to be connected to $\RR\SS$ is randomly selected. In case the resulting solution is infeasible, the procedure continues by selecting a new random vertex of $\RR\SS$ as insertion point. 
\end{itemize}

\section{Computational results} \label{sec:Results} 

In this section, we provide extensive computational tests to evaluate the effectiveness of the proposed ALNS. All algorithms have been implemented in C++ and tested on a PC with a 2.8 GHz Intel Dual Core and 4 GB of RAM. We used ILOG Cplex 12.3 to run the RCKRMSF model of Section \ref{sec:Modeling}. We present results on new RCKRMSF instances that we randomly generated and on DCMST benchmark instances from the literature. The ALNS has been executed five times on each attempted instance. Due to the large number of tests, in this section we only report aggregated results, so as to facilitate comprehension. 

\subsection{Instance generation} \label{sec:Data} 

We generated two sets of instances, the first having Euclidean distances and the second non-Euclidean. The first set contains 180 instances having number of vertices $n \in \{20, 30, 40, 50\}$ and number of sources $m\in \left\{2,3,4\right\}$. Three values were selected for the minimum path reliability, namely $\alpha \in \left\{0.95, 0.90, 0.80\right\}$, and the reliabilities of the arcs were created as random numbers in the interval $\left[\alpha, 1\right]$. The coordinates of the vertices were randomly selected in the square $\left[0,100\right]\times\left[0,100\right]$, and, for each $i,j \in V$, $c_{ij}$ was set equal to the Euclidean distance between the coordinates of $i$ and $j$. To take randomness into account, for each combination of $n$, $m$ and $\alpha$, five random instances have been generated, differing one another in the distances and reliabilities on the arcs.
The second set was created mainly to compare Euclidean with non-Euclidean distances. This set was indeed generated as the first one for what concerns $n$, $m$, $\alpha$ and arc reliabilities, but it contains distances that have been randomly generated for each arc in the interval $\left[1,100\right]$. For each combination of $n$, $m$ and $\alpha$, two random instances have been generated. 

All instances that we generated are provided as on-line supplemental material to this paper.
 
\subsection{Parameter tuning} \label{sec:Parameter} 

The ALNS of Section \ref{sec:Solutionmethod}  contains a number of parameters that have to be finely tuned. To find their best configuration, we opted to use a simple {\em genetic algorithm} (GA). GAs (see, e.g., \citealt{SGK14}) are well suited for the optimization of complex problems, and parameter tuning is one of the fields in which they have been successfully applied, as recently shown in, e.g., \cite{LNWF14, YX14}. The sketch of the GA that we adopted for selecting the best ALNS parameters is described in Algotithm \ref{sec:Alg3}.

\begin{table}[htb]
\begin{algorithm}[H]
\hrulefill \\
\footnotesize
\label{sec:Alg3}
\begin{enumerate}
	\item Generate an initial random population of individuals (each being a combination of parameters)
	\item Calculate the fitness of each individual
	\item Produce the next population by using the following operators
	\begin{itemize}
		\item Selection: Individuals with high fitness values are selected to breed a new generation
		\item Crossover: Groups of two parents are selected and combined to create offspring solutions 
		\item Mutation: Each solution is locally randomly modified to increase diversity	
	\end{itemize}
	\item Repeat steps 2 and 3 until the stopping criterion is met
	\item Return the last population
\end{enumerate}
\hrulefill
\caption{\small The general framework of the GA adopted for the ALNS parameter tuning}
\end{algorithm}
\end{table}

For this preliminary tuning, a test set of 20 instances having 20 or 40 vertices have been randomly generated as in the first set of Section \ref{sec:Data}. The initial GA population was set to 20 random individuals. The chromosome of each individual consists of a list of the ALNS parameters to be tuned, and its fitness is the average value of the 20 solutions found by executing the ALNS with this list of parameters on the 20 instances in the test set.
We adopted a simple truncation selection, in which the top $(1/s)^{th}$ of the individuals get $s$ copies each in the mating pool, with $s=2$ in our tests. One-point crossover was chosen as the operator for recombining the two parent chromosomes (a point is randomly selected  over the chromosome length, and all data beyond this point is swapped between the two parents). As mutation procedure, we opted to change an arbitrary bit in a chromosome from its original state. We decided to stop the GA after the generation of 30 populations. 

Table 1 shows the list of ALNS parameters, their brief descriptions, the sets of parameter values that have been tested with the GA, and, in the last column, the final values that have been adopted on the basis of the results. It is worth mentioning that the performance of the ALNS was not very sensitive to the set of parameters. We applied indeed the Wilcoxon signed rank test, a well-known nonparametric statistical test from \cite{W45}, to perform a pairwise comparison of parameters. To do so, we selected the two parameter combinations having, respectively, the best and the worst average performance (by considering the average objective function over the 20 test instances). The Wilcoxon test showed that these combinations do not dominate each other at a significance level of $\alpha < 0.05$. 
This confirms the robustness of the developed ALNS. 

\begin{table}
\begin{center}
\tabcolsep=3.5pt
\scriptsize
\label{tab:table1}
\caption{\small ALNS parameters}  
\begin{tabular}{llll} 
\toprule 
       Parameter            & Description       																							 & Different tested values       & Adopted value   \\ 
 \cmidrule(l){1-1}\cmidrule(l){2-4} 
       $\theta$              &initial temperature 													 & \{500, 1000, 2000, 5000, 10000, 15000\} 	   & 1000         \\
       $\epsilon$            &final temperature   										       & \{0.0001\}     				      & 0.0001       \\ 
       $\lambda$             &cooling rate       													 & \{0.1, 0.01, 0.001, 0.0001, 0.00001\}        & 0.1          \\
       $\eta$                &reaction factor    													 & \{0.1, 0.2, 0.3, 0.4, 0.5\}                	& 0.2          \\ 
       $\CC$                 &distance between two vertices in shaking 3 & \{20, 30, 40, 50, 60\}  & 50           \\ 
       $\phi_1$              &number of segments          							   & \{3, 5, 10, 20, 30, 40\}          		      & 5            \\
       $\phi_2$              &number of iterations in each segment         & \{3, 5, 10, 20, 30, 40\}  		        	    & 10           \\
			 $\zeta_1$             & $\pi$ increment value 1 & \{30, 50, 70, 90, 110, 130\} & 50\\
			 $\zeta_2$             & $\pi$ increment value 2 & \{10, 20, 30, 40, 50, 60\}& 20\\
			 $\zeta_3$             & $\pi$ increment value 3 & \{3, 5, 10, 15, 20, 25\}& 5\\
\bottomrule
\end{tabular}
\end{center}
\end{table}

\subsection{Computational results on DCMST benchmarks} \label{sec:Results3} 

With the aim of assessing the effectiveness of the  ALNS, we considered the DCMST discussed in Section \ref{sec:review}. The DCMST is a special case of the RCKRMSF in which a single source is used to serve demand vertices and a maximum delay is allowed on each path. Our ALNS can be easily adapted to solve this problem by considering delays instead of reliabilities on the arcs. 

A relevant work on the DCMST is that of \cite{GPS08}, who proposed a number of instances and solved most of them to proven optimality in one day of CPU time by means of algorithms based on Lagrangian relaxation and column generation. The authors considered instances with 20 and 40 demand vertices, divided into three main groups, namely TR, TC and TE. While the TR group contains instances for which the edge costs are random values, the TC and TE groups both include instances with Euclidean costs, with the only difference that the source vertex is near the center of the grid of vertices for the TC group and near the border for the TE.  Instances are denoted as X,Y [Z,W], where X is the group, Y is the number of demand vertices and [Z,W] gives the range in which the arc delays have been created, with Z=1 and W $\in \{2, 5, 10, 100, 1000\}$. For each such configuration, four instances have been generated by using different values on the maximum delay bound, resulting in a total of 120 instances. The ALNS was executed five times on each instance.

Similarly to \cite{GPS08}, we present separate results for the classes involving small and large delays. The former results are in Table \ref{tab:table9-A} and the latter ones in Table \ref{tab:table9-B}. For each group of four instances per line, each table provides the average (over the four instances per line) of the best (over the five attempts) gap between the ALNS solution value and the best known solution value reported by \cite{GPS08}, and the average ALNS execution time (over the four instances and five attempts). The ALNS was executed with the default parameters proposed for the RCKRMSF.  It is clear from the tables that the ALNS is very effective in tackling the DCMST, being able to achieve very low gaps in a matter of seconds. The average best gap is usually below 0.5\% and just in two cases (both for TR,40) reaches 0.71\%. Essentially, the average gap over all the available DCMST instances, with respect to the best know solutions, is around 0.21\%, and the average CPU time is just 2.38 seconds.

\begin{table}[htb]
\begin{center}
\scriptsize
\caption{\small ALNS comparison with best solutions on DCMST benchmarks -- small delays (4 instances per line)} \label{tab:table9-A}
\begin{tabular}{lrrlrrlrr} 
\toprule
Class 1 & Avg. & Avg. & Class 2 & Avg.  & Avg. & Class 3 & Avg.  & Avg.   \\ 
  & best gap & time &  & best gap &  time &  & best gap & time  \\ 
 \cmidrule(l){1-3}\cmidrule(l){4-6}\cmidrule(l){7-9}   
 TR,20 [1,2] & 0.00 & 1.61 & TR,20 [1,5] & 0.00 & 1.56 & TR,20 [1,10] & 0.00 & 1.64 \\  
 TC,20 [1,2] & 0.00 & 1.59 & TC,20 [1,5] & 0.00 & 1.68 & TC,20 [1,10] & 0.08 & 1.63  \\  
 TE,20 [1,2] & 0.00 & 1.27 & TE,20 [1,5] & 0.00 & 1.30 & TE,20 [1,10] & 0.03 & 1.31  \\  
 TR,40 [1,2] & 0.40 & 2.82 & TR,40 [1,5] & 0.52 & 3.16 & TR,40 [1,10] & 0.47 & 3.06  \\  
 TC,40 [1,2] & 0.11 & 2.82 & TC,40 [1,5] & 0.22 & 2.98 & TC,40 [1,10] & 0.22 & 3.08 \\  
 TE,40 [1,2] & 0.25 & 3.25 & TE,40 [1,5] & 0.64 & 3.24 & TE,40 [1,10] & 0.35 & 3.32  \\ 
\bottomrule
\end{tabular}
\end{center}
\end{table}

\begin{table}[htb]
\begin{center}
\tabcolsep=8pt
\scriptsize
\caption{\small ALNS comparison with best solutions on DCMST benchmarks -- large delays (4 instances per line)} \label{tab:table9-B}
\begin{tabular}{lrrlrr} 
\toprule
Class 4 & Avg. & Avg. & Class 5 & Avg.  & Avg.    \\ 
  & best gap & time &  & best gap &  time   \\ 
 \cmidrule(l){1-3}\cmidrule(l){4-6}
 TR,20 [1,100] & 0.35 & 1.69 & TR,40 [1,1000] & 0.20 & 1.76 \\  
 TC,20 [1,100] & 0.06 & 1.81 & TC,40 [1,1000] & 0.00 & 1.86 \\  
 TE,20 [1,100] & 0.00 & 1.29 & TE,40 [1,1000] & 0.00 & 1.35 \\  
 TR,40 [1,100] & 0.71 & 3.44 & TR,40 [1,1000] & 0.71 & 2.57 \\  
 TC,40 [1,100] & 0.17 & 3.05 & TC,40 [1,1000] & 0.25 & 2.19 \\  
 TE,40 [1,100] & 0.74 & 3.62 & TE,40 [1,1000] & 0.52 & 3.72 \\ 
\bottomrule
\end{tabular}
\end{center}
\end{table}

\begin{figure}[h!]
  \caption{The effect of parameter $\phi_2$ on the average best gap for the three DCMST groups}
  \centering
    \includegraphics[angle=0,width=110mm]{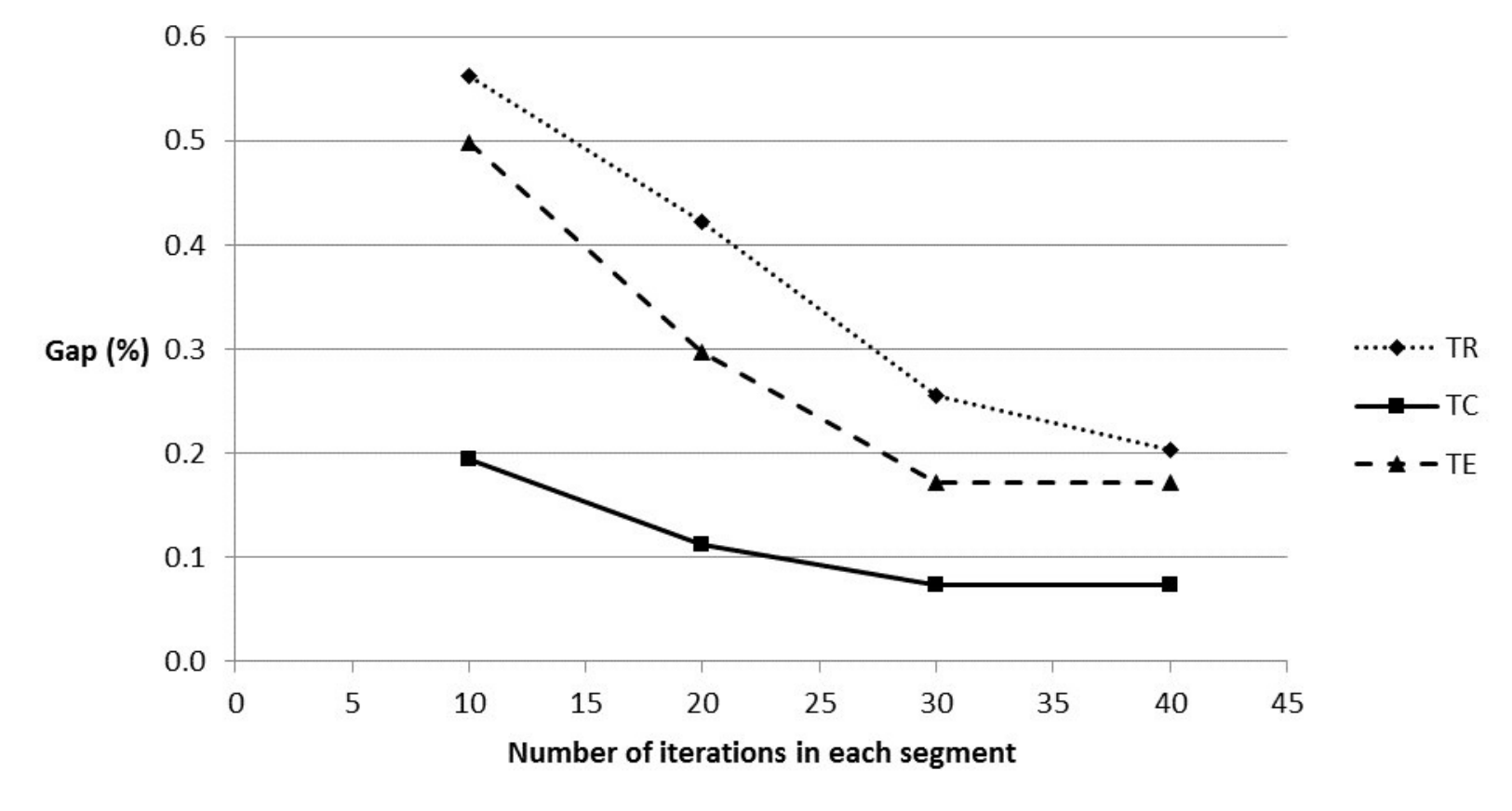}
    \label{fig:SAI}
\end{figure}

\begin{figure}[h!]
  \caption{The effect of parameter $\phi_1$ on the average best gap or the three DCMST groups}
  \centering
    \includegraphics[angle=0,width=110mm]{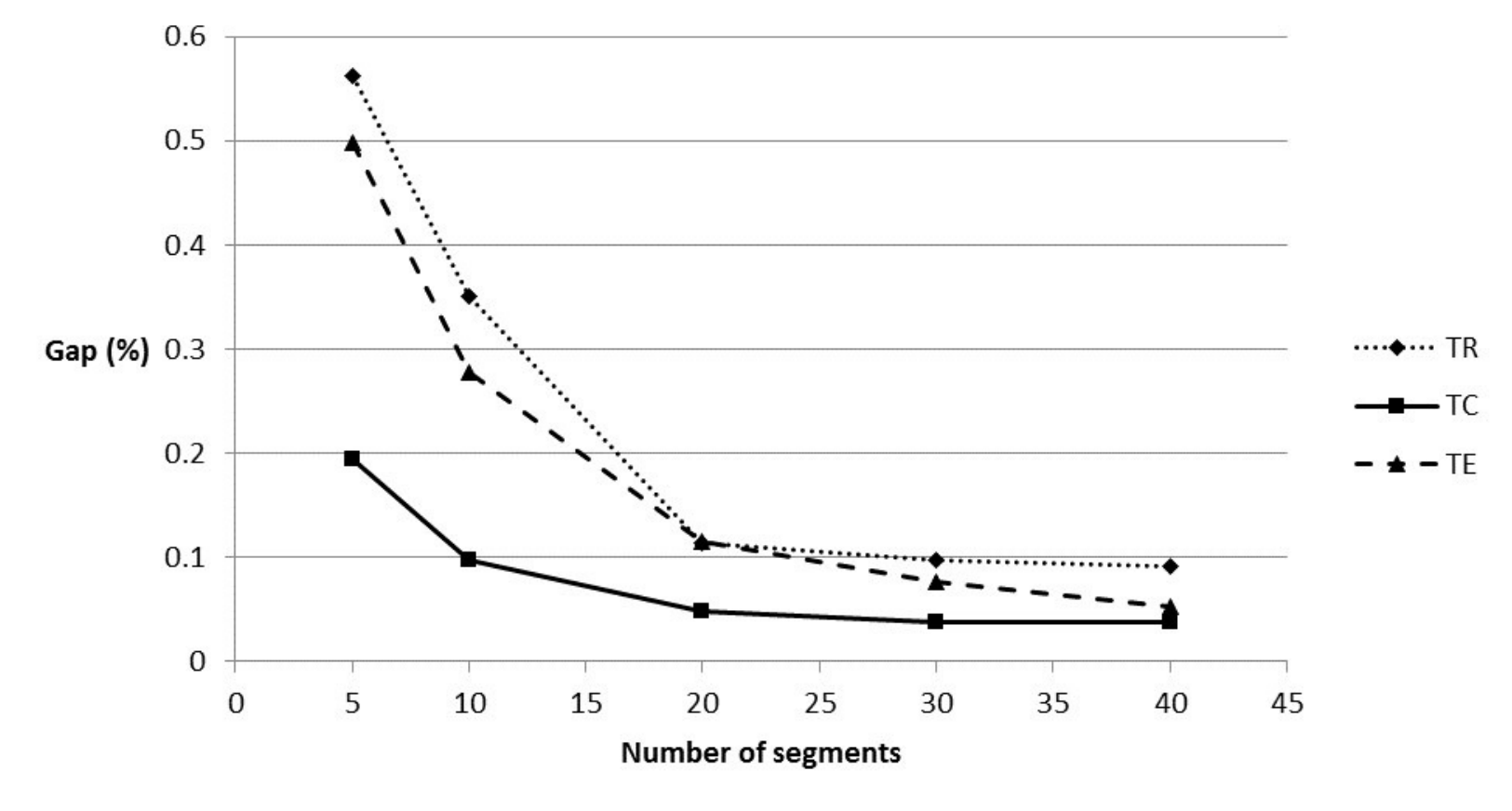}
    \label{fig:SAII}
\end{figure}

We use this set of instances also to provide a sensitivity analysis for $\phi_1$ (number of global ALNS iterations) and $\phi_2$ (number of inner iterations). Figure \ref{fig:SAI} shows the average best gap for all DCMST instances of the three groups with respect to different values of $\phi_2$, by keeping $\phi_1 = 5$. It can be noticed that the gap tends to zero by increasing the value of $\phi_2$, showing that the ALNS tends not to be stuck in local optima. Figure \ref{fig:SAII} 
shows the same analysis for the case in which $\phi_2 = 10$ and $\phi_1$ takes different values ranging from 5 to 40. Also in this case, by increasing the number of iterations the optimality gap for all TR, TC and TE groups decreases and gets below 0.1\%. 

\subsection{Computational results on Euclidean RCKRMSF instances} \label{sec:Results1} 
In this section, we focus on the solution of the new RCKRMSF instances that we generated. Tables \ref{tab:table3}, \ref{tab:table4} and \ref{tab:table5} give the results that we obtained on Euclidean instances with $\alpha$=0.95, 0.90 and 0.80, respectively, by running the MILP model composed of \eqref{eq:obj1}--\eqref{eq:con10} and \eqref{eq:con12}--\eqref{eq:con15} with CPLEX and the  ALNS of Section \ref{sec:Solutionmethod}. CPLEX was allowed a maximum run time of 18000 CPU seconds. Each line in the tables report average results on the selected group of five instances having the same values of $n$, $m$ and $\alpha$. For the MILP model, we provide the lower bound value reported at the end of the run in column ``LB'', the solution value in column ``UB'', and the elapsed CPU seconds in column ``time''. For the ALNS, which was run five times on each instance, we provide, in order, the best, average and worst solution values, the average and best percentage gaps from the UB of the MILP model, and the average execution time. The gaps have been computed as 100$\times$(UB$_{ALNS}$$-$UB$_{MILP}$)/UB$_{MILP}$.

\begin{table}[h!]
\begin{center}
\tabcolsep=8pt
\scriptsize
\caption{\small Computational results on Euclidean instances with $\alpha=0.95$ (5 instances per line)}  \label{tab:table3}
\begin{tabular}{rrrrrrrrrrr} 
\toprule 
 \textit{n} & \textit{m} & \multicolumn{3}{c}{MILP model} & \multicolumn{6}{c}{{ALNS}} \\
 \cmidrule(l){1-2}\cmidrule(l){3-5}\cmidrule(l){6-11}   
 &  & LB & UB & Time & Best & Worst & Avg. & Avg. & Best & Avg. \\
 &  &    &       &      &  UB    &   UB    &  UB       & gap  & gap  & time \\
 \cmidrule(l){1-2}\cmidrule(l){3-5}\cmidrule(l){6-11}   
20 & 2 & 325.82 & 325.82 & 531.83 & 325.82 & 329.83 & 326.62 & 0.23 & 0.00 & 1.40 \\   
20 & 3 & 281.06 & 281.06 & 201.60 & 281.06 & 281.06 & 281.06 & 0.00 & 0.00 & 1.72 \\   
20 & 4 & 255.60 & 255.60 & 106.52 & 255.60 & 256.05 & 255.96 & 0.13 & 0.00 & 1.36 \\  
30 & 2 & 418.77 & 445.85 & 18000.00 & 446.04 & 447.46 & 446.58 & 0.16 & 0.04 & 1.76 \\  
30 & 3 & 394.98 & 420.63 & 18000.00 & 421.42 & 423.89 & 422.71 & 0.46 & 0.17 & 1.44 \\  
30 & 4 & 366.06 & 383.04 & 15701.08 & 383.04 & 387.12 & 384.71 & 0.45 & 0.00 & 1.69  \\  
40 & 2 & 483.48 & 556.20 & 18000.00 & 527.58 & 539.12 & 533.44 & $-$3.86 & $-$4.93 & 2.57 \\  
40 & 3 & 447.38 & 490.20 & 18000.00 & 488.15 & 494.41 & 490.49 & 0.08 & $-$0.40 & 2.24  \\  
40 & 4 & 425.12 & 460.64 & 18000.00 & 459.61 & 463.28 & 460.64 & 0.01 & $-$0.20 & 2.45  \\  
50 & 2 & 524.96 & 677.13 & 18000.00 & 605.31 & 617.26 & 610.20 & $-$9.52 & $-$10.23 & 5.88  \\  
50 & 3 & 495.23 & 622.93 & 18000.00 & 567.68 & 576.16 & 570.42 & $-$8.24 & $-$8.67 & 4.56 \\  
50 & 4 & 466.33 & 567.37 & 18000.00 & 530.30 & 538.87 & 533.77 & $-$5.77 & $-$6.37 & 4.68  \\  
 \cmidrule(l){1-2}\cmidrule(l){3-5}\cmidrule(l){6-11}   
\multicolumn{2}{c}{Average} & 407.07 & 457.20 & 13378.42 & 440.97 & 446.21 & 443.05 & $-$2.16 & $-$2.55 & 2.65 \\  
\bottomrule
\end{tabular}
\end{center}
\end{table}
\begin{table}[h!]
\begin{center}
\tabcolsep=8pt
\scriptsize
\caption{\small Computational results on Euclidean instances with $\alpha=0.90$ (5 instances per line)} \label{tab:table4}
\begin{tabular}{rrrrrrrrrrr} 
\toprule 
 \textit{n} & \textit{m} & \multicolumn{3}{c}{MILP model} & \multicolumn{6}{c}{{ALNS}} \\
 \cmidrule(l){1-2}\cmidrule(l){3-5}\cmidrule(l){6-11}   
 &  & LB & UB & Time & Best & Worst & Avg. & Avg. & Best & Avg. \\
 &  &    &       &      &  UB    &   UB    &  UB       & gap  & gap  & time \\
 \cmidrule(l){1-2}\cmidrule(l){3-5}\cmidrule(l){6-11}   
20 & 2 & 332.77 & 332.77 & 431.80 & 332.77 & 334.93 & 333.21 & 0.12 & 0.00 & 1.14 \\   
20 & 3 & 293.89 & 293.89 & 166.57 & 293.89 & 294.65 & 294.19 & 0.10 & 0.00 & 1.30 \\   
20 & 4 & 266.55 & 266.55 & 44.40 & 266.55 & 267.64 & 266.77 & 0.08 & 0.00 & 1.44 \\  
30 & 2 & 402.61 & 444.93 & 18000.00 & 443.75 & 453.15 & 446.56 & 0.34 & $-$0.27 & 1.62 \\  
30 & 3 & 382.78 & 406.83 & 16023.94 & 406.83 & 408.20 & 407.51 & 0.17 & 0.00 & 1.65 \\  
30 & 4 & 358.15 & 376.43 & 14631.73 & 376.43 & 377.71 & 376.85 & 0.11 & 0.00 & 1.67  \\  
40 & 2 & 483.34 & 551.44 & 18000.00 & 537.47 & 543.89 & 540.42 & $-$1.97 & $-$2.52 & 2.61 \\  
40 & 3 & 451.25 & 497.56 & 18000.00 & 496.57 & 501.25 & 498.68 & 0.22 & $-$0.20 & 2.62  \\  
40 & 4 & 429.33 & 474.45 & 18000.00 & 472.62 & 477.42 & 474.29 & $-$0.03 & $-$0.37 & 2.59  \\  
50 & 2 & 524.00 & 702.91 & 18000.00 & 585.99 & 602.61 & 593.08 & $-$14.70 & $-$15.65 & 5.38  \\  
50 & 3 & 495.10 & 593.27 & 18000.00 & 550.57 & 560.99 & 554.76 & $-$6.14 & $-$6.82 & 4.92 \\  
50 & 4 & 466.72 & 560.15 & 18000.00 & 517.15 & 523.04 & 519.52 & $-$6.90 & $-$7.31 & 4.55  \\  
 \cmidrule(l){1-2}\cmidrule(l){3-5}\cmidrule(l){6-11}   
\multicolumn{2}{c}{Average} & 407.21 & 458.43 & 13108.20 & 440.05 & 445.46 & 442.15 & $-$2.38 & $-$2.76 & 2.62 \\    
\bottomrule
\end{tabular}
\end{center}
\end{table}
\begin{table}[h!]
\begin{center}
\scriptsize
\caption{\small Computational results on Euclidean instances with $\alpha=0.80$ (5 instances per line)}  \label{tab:table5}
\begin{tabular}{rrrrrrrrrrr} 
\toprule 
 \textit{n} & \textit{m} & \multicolumn{3}{c}{MILP model} & \multicolumn{6}{c}{{ALNS}} \\
 \cmidrule(l){1-2}\cmidrule(l){3-5}\cmidrule(l){6-11}   
 &  & LB & UB & Time & Best & Worst & Avg. & Avg. & Best & Avg. \\
 &  &    &       &      &  UB    &   UB    &  UB       & gap  & gap  & time \\
 \cmidrule(l){1-2}\cmidrule(l){3-5}\cmidrule(l){6-11}   
20 & 2 & 323.55 & 323.55 & 521.29 & 323.55 & 324.82 & 324.06 & 0.15 & 0.00 & 1.38 \\   
20 & 3 & 289.07 & 289.07 & 599.96 & 289.07 & 290.41 & 290.02 & 0.27 & 0.00 & 1.42 \\   
20 & 4 & 264.06 & 264.06 & 409.47 & 264.06 & 264.79 & 264.55 & 0.16 & 0.00 & 1.51 \\  
30 & 2 & 412.56 & 455.50 & 18000.00 & 455.33 & 461.63 & 457.85 & 0.51 & $-$0.02 & 1.89 \\  
30 & 3 & 390.60 & 418.67 & 18000.00 & 419.20 & 419.98 & 419.51 & 0.20 & 0.13 & 1.84 \\  
30 & 4 & 363.04 & 379.28 & 16189.15 & 379.28 & 379.96 & 379.68 & 0.11 & 0.00 & 1.73  \\  
40 & 2 & 489.87 & 559.68 & 18000.00 & 546.71 & 559.95 & 552.10 & $-$3.37 & $-$4.30 & 2.60 \\  
40 & 3 & 453.61 & 503.92 & 18000.00 & 499.72 & 507.78 & 502.80 & $-$0.15 & $-$0.74 & 2.92  \\  
40 & 4 & 429.22 & 479.47 & 18000.00 & 474.93 & 479.34 & 476.88 & $-$0.47 & $-$0.90 & 2.61  \\  
50 & 2 & 514.04 & 647.95 & 18000.00 & 584.71 & 596.33 & 590.69 & $-$8.56 & $-$9.47 & 5.20  \\  
50 & 3 & 484.83 & 593.58 & 18000.00 & 548.41 & 562.12 & 552.54 & $-$6.82 & $-$7.49 & 4.70 \\  
50 & 4 & 455.51 & 563.24 & 18000.00 & 516.11 & 517.95 & 516.97 & $-$7.78 & $-$7.93 & 4.81  \\  
\cmidrule(l){1-2}\cmidrule(l){3-5}\cmidrule(l){6-11}   
\multicolumn{2}{c}{Average} & 405.83 & 456.50 & 13476.66 & 441.76 & 447.09 & 443.97 & $-$2.14 & $-$2.56 & 2.72 \\    
\bottomrule
\end{tabular}
\end{center}
\end{table}

Out of the 180 instances addressed in the three tables, {49} instances were solved to proven optimality by the MILP model in the given time limit. These include all 45 instances with 20 demand vertices, and 4 with 30 vertices. For the unsolved instances, the difference between the UB and LB values produced by the model can be quite high. On average, the overall gap between these two values is slightly above 10\%, for all values of $\alpha$.

The ALNS finds all 49 proven optimal UBs in at least one of the five runs performed on those instances. On the other instances, it usually produces very good gaps. The best gap varies from  0.17\% (for $n=30$, $m=3$ and $\alpha = 0.95$) to $-$15.65\% (for $n=50$, $m=2$ and $\alpha = 0.90$). On average, the best gap shows improvement  of about 2.5\% with respect to the UBs produced by the MILP model, but such improvements can be very relevant on instances with 50 vertices, with values ranging between $-$6\% and $-$15\%. The good performance of the ALNS can also be noticed by looking at the worst and average UBs, which are not far away from the best ones. The overall average gap is between $-$2.1\% and $-$2.4\%, and is quite consistent for the different values taken by $\alpha$. These good-quality solutions are obtained in very short times, involving less than 3 seconds on average and nearly 6 seconds in the worst case (for $n=50$, $m=2$ and $\alpha = 0.95$). By looking at these tables, we can conclude that the MILP model is effective up to 20 demand vertices, whereas the ALNS is effective and fast on all instances. It can also be noted that the algorithms are not particularly sensitive to the value of $\alpha$.

Tables \ref{tab:table6}, \ref{tab:table7} and \ref{tab:table8} provide average results on the non-Euclidean set of instances, focusing again on $\alpha$=0.95, 0.90 and 0.80, respectively. The values reported have the same meanings as those in  Tables \ref{tab:table3}--\ref{tab:table5}, but now refer to just two instances per line instead of five.
Out of 72 instances, the ALNS outperforms the UB value obtained by CPLEX in 36 cases at least once in the five runs for instance. For the remaining instances, it always finds the same objective value of the MILP model, still at least once in the five runs. Overall, all best gaps are always null or negative, whereas the average gaps might rarely take positive values, although these are very small. More variability with respect to the results on the Euclidean instances can be noticed on the overall average gaps, which vary from a minimum of $-$5.27\% for $\alpha = 0.90$ to a maximum of $-$2.22\% for $\alpha = 0.80$. This might be induced by the reduced number of instances in the second test bed. The MILP model is now a bit more effective on the instances having 30 demand vertices, but remains very slow for the larger instances. Also for this test bed, we can conclude that the ALNS is fast and effective.
\begin{table}[h!]
\begin{center}
\scriptsize
\caption{\small Computational results on non-Euclidean instances with $\alpha=0.95$ (2 instances per line)}  \label{tab:table6}
\begin{tabular}{rrrrrrrrrrr} 
\toprule 
 \textit{n} & \textit{m} & \multicolumn{3}{c}{Mathematical model} & \multicolumn{6}{c}{{ALNS}} \\
 \cmidrule(l){1-2}\cmidrule(l){3-5}\cmidrule(l){6-11}   
 &  & LB & UB & Time & Best & Worst & Avg. & Avg. & Best & Avg. \\
 &  &    &       &      &  UB    &   UB    &  UB       & gap  & gap  & time \\
 \cmidrule(l){1-2}\cmidrule(l){3-5}\cmidrule(l){6-11}   
20 & 2 & 162.50 & 162.50 & 93.31 & 162.50 & 164.00 & 163.20 & 0.35 & 0.00 & 0.82 \\   
20 & 3 & 144.00 & 144.00 & 134.52 & 144.00 & 144.00 & 144.00 & 0.00 & 0.00 & 0.89 \\   
20 & 4 & 132.50 & 132.50 & 64.05 & 132.50 & 133.00 & 132.60 & 0.06 & 0.00 & 0.86 \\  
30 & 2 & 199.91 & 208.00 & 14578.21 & 207.50 & 208.00 & 207.60 & $-$0.16 & $-$0.20 & 2.05 \\  
30 & 3 & 176.50 & 176.50 & 8939.29 & 176.50 & 176.50 & 176.50 & 0.00 & 0.00 & 2.10 \\  
30 & 4 & 163.00 & 163.00 & 4301.67 & 163.00 & 163.00 & 163.00 & 0.00 & 0.00 & 1.73  \\  
40 & 2 & 177.14 & 219.00 & 18000.00 & 212.00 & 213.50 & 212.30 & $-$3.06 & $-$3.20 & 2.36 \\  
40 & 3 & 167.55 & 207.00 & 18000.00 & 199.50 & 202.50 & 200.30 & $-$3.24 & $-$3.63 & 2.32  \\  
40 & 4 & 157.36 & 195.00 & 18000.00 & 186.50 & 190.50 & 187.90 & $-$3.34 & $-$4.10 & 2.06  \\  
50 & 2 & 160.72 & 239.00 & 18000.00 & 204.00 & 210.50 & 207.10 & $-$13.05 & $-$14.29 & 4.42  \\  
50 & 3 & 151.06 & 205.50 & 18000.00 & 189.50 & 196.00 & 192.30 & $-$5.80 & $-$7.22 & 3.76 \\  
50 & 4 & 144.37 & 204.50 & 18000.00 & 180.50 & 190.50 & 183.70 & $-$8.88 & $-$10.62 & 3.64  \\  
 \cmidrule(l){1-2}\cmidrule(l){3-5}\cmidrule(l){6-11}   
\multicolumn{2}{c}{Average} & 161.38 & 188.04 & 11342.59 & 179.83 & 182.67 & 180.88 & $-$3.09 & $-$3.60 & 2.25 \\    
\bottomrule
\end{tabular}
\end{center}
\end{table}
\begin{table}[h!]
\begin{center}
\scriptsize
\caption{\small Computational results on non-Euclidean instances with $\alpha=0.90$ (2 instances per line)}  \label{tab:table7}
\begin{tabular}{rrrrrrrrrrr} 
\toprule 
 \textit{n} & \textit{m} & \multicolumn{3}{c}{Mathematical model} & \multicolumn{6}{c}{{ALNS}} \\
 \cmidrule(l){1-2}\cmidrule(l){3-5}\cmidrule(l){6-11}   
 &  & LB & UB & Time & Best & Worst & Avg. & Avg. & Best & Avg. \\
 &  &    &       &      &  UB    &   UB    &  UB       & gap  & gap  & time \\
 \cmidrule(l){1-2}\cmidrule(l){3-5}\cmidrule(l){6-11}   
20 & 2 & 148.50 & 148.50 & 17.96 & 148.50 & 148.50 & 148.50 & 0.00 & 0.00 & 1.02 \\   
20 & 3 & 128.00 & 128.00 & 13.90 & 128.00 & 128.00 & 128.00 & 0.00 & 0.00 & 0.95 \\   
20 & 4 & 118.00 & 118.00 & 13.38 & 118.00 & 118.00 & 118.00 & 0.00 & 0.00 & 1.39 \\  
30 & 2 & 203.31 & 212.50 & 17624.33 & 211.50 & 217.50 & 213.90 & 0.79 & $-$0.44 & 1.61 \\  
30 & 3 & 185.27 & 196.50 & 18000.00 & 196.50 & 198.00 & 196.80 & 0.14 & 0.00 & 2.61 \\  
30 & 4 & 170.89 & 183.50 & 14578.38 & 183.50 & 183.50 & 183.50 & 0.00 & 0.00 & 2.19  \\  
40 & 2 & 204.46 & 254.00 & 18000.00 & 245.50 & 249.00 & 248.30 & $-$2.24 & $-$3.31 & 2.68 \\  
40 & 3 & 188.65 & 252.00 & 18000.00 & 231.00 & 232.50 & 231.30 & $-$7.95 & $-$8.06 & 2.92  \\  
40 & 4 & 175.26 & 207.50 & 18000.00 & 205.00 & 205.00 & 205.00 & $-$1.11 & $-$1.11 & 2.17  \\  
50 & 2 & 157.98 & 299.50 & 18000.00 & 193.00 & 198.00 & 195.30 & $-$34.82 & $-$35.59 & 5.38  \\  
50 & 3 & 151.16 & 218.00 & 18000.00 & 184.50 & 190.50 & 187.30 & $-$14.11 & $-$15.38 & 3.94 \\  
50 & 4 & 142.21 & 184.00 & 18000.00 & 175.00 & 179.50 & 176.80 & $-$3.93 & $-$4.90 & 3.97  \\  
 \cmidrule(l){1-2}\cmidrule(l){3-5}\cmidrule(l){6-11}   
\multicolumn{2}{c}{Average} & 164.47 & 200.17 & 13187.33 & 185.00 & 187.33 & 186.06 & $-$5.27 & $-$5.73 & 2.57 \\      
\bottomrule
\end{tabular}
\end{center}
\end{table}
\begin{table}[h!]
\begin{center}
\tabcolsep=8pt
\scriptsize
\caption{\small Computational results on non-Euclidean instances with $\alpha=0.80$ (2 instances per line)}  \label{tab:table8}
\begin{tabular}{rrrrrrrrrrr} 
\toprule 
 \textit{n} & \textit{m} & \multicolumn{3}{c}{Mathematical model} & \multicolumn{6}{c}{{ALNS}} \\
 \cmidrule(l){1-2}\cmidrule(l){3-5}\cmidrule(l){6-11}   
 &  & LB & UB & Time & Best & Worst & Avg. & Avg. & Best & Avg. \\
 &  &    &       &      &  UB    &   UB    &  UB       & gap  & gap  & time \\
 \cmidrule(l){1-2}\cmidrule(l){3-5}\cmidrule(l){6-11} 
20 & 2 & 186.00 & 186.00 & 73.26 & 186.00 & 186.00 & 186.00 & 0.00 & 0.00 & 0.84 \\   
20 & 3 & 145.00 & 145.00 & 103.29 & 145.00 & 145.00 & 145.00 & 0.00 & 0.00 & 1.02 \\   
20 & 4 & 117.00 & 117.00 & 6.54 & 117.00 & 117.00 & 117.00 & 0.00 & 0.00 & 1.33 \\  
30 & 2 & 189.89 & 204.00 & 10162.24 & 204.00 & 209.50 & 205.80 & 0.98 & 0.00 & 2.07 \\  
30 & 3 & 169.50 & 169.50 & 6125.59 & 169.50 & 170.00 & 169.90 & 0.31 & 0.00 & 1.72 \\  
30 & 4 & 158.00 & 158.00 & 6984.14 & 158.00 & 158.50 & 158.30 & 0.24 & 0.00 & 1.88  \\  
40 & 2 & 195.31 & 238.50 & 18000.00 & 231.50 & 234.50 & 232.50 & $-$2.66 & $-$3.07 & 2.47 \\  
40 & 3 & 181.05 & 244.00 & 18000.00 & 224.00 & 225.50 & 224.60 & $-$7.86 & $-$8.11 & 2.20  \\  
40 & 4 & 170.96 & 211.50 & 18000.00 & 203.00 & 204.00 & 203.30 & $-$3.89 & $-$4.02 & 2.38  \\  
50 & 2 & 167.49 & 230.00 & 18000.00 & 212.50 & 217.50 & 213.90 & $-$7.01 & $-$7.61 & 4.97  \\  
50 & 3 & 158.77 & 206.50 & 18000.00 & 197.00 & 207.00 & 199.10 & $-$3.32 & $-$4.28 & 3.87 \\  
50 & 4 & 148.10 & 191.50 & 18000.00 & 183.50 & 189.50 & 184.80 & $-$3.45 & $-$4.14 & 3.96  \\  
 \cmidrule(l){1-2}\cmidrule(l){3-5}\cmidrule(l){6-11}   
\multicolumn{2}{c}{Average} & 165.59 & 191.79 & 10954.59 & 185.92 & 188.67 & 186.68 & $-$2.22 & $-$2.60 & 2.39 \\      
\bottomrule
\end{tabular}
\end{center}
\end{table}

\section{Conclusions} \label{sec:Conclusion} 
We proposed the reliability constraint \textit{k}-rooted minimum spanning forest. The goal of the problem is to find a \textit{k}-rooted minimum cost forest that satisfies a minimum required reliability on each path from a source to a demand vertex. The problem has relevant applications in the design of several types of networks. We solved it by means of a mathematical model and an adaptive large neighborhood search metaheuristic. The model is effective in solving instances with up to 20 demand vertices, whereas the metaheuristic is both fast and effective not only on all instances of the proposed problem, but also on the special case known in the literature as the delay constrained minimum spanning tree. Such conclusions are motivated by extensive computational tests on randomly created instances that have been now made publicly available to stimulate further research. 

Instances with just 30 demand vertices remain unsolved to proven optimality despite relevant computing efforts, so, as future research, much could be done in terms of exact algorithms to attempt solving larger instances. Another interesting research direction is that of studying the budget-constrained version of the problem, that is, the problem variant in which the goal is to construct a maximum reliability spanning forest with the presence of a limitation on the available network construction cost. 

\singlespacing
\bibliography{biblio}
\bibliographystyle{mmsbib}

\end{document}